\input amstex
\documentstyle{amsppt}
\magnification=1200
\vsize 9.0truein
\hsize 6.5truein
%\NoBlackBoxes
\leftheadtext{Montgomery \& Soundararajan}
\loadeusm
\loadbold
\topmatter
\title
Primes in short intervals
\endtitle
\author
Hugh L. Montgomery\footnotemark \\
K. Soundararajan\footnotemark
\endauthor
\dedicatory
Dedicated to Freeman Dyson, with best wishes on the occasion of his eightieth
birthday.
\enddedicatory
\abstract
Contrary to what would be predicted on the basis of Cram\'er's model concerning
the distribution of prime numbers, we develop evidence that the distribution of
$\psi(x+H)- \psi(x)$, for $0\le x\le N$, is approximately normal with mean $\sim H$
and variance $\sim H\log N/H$, when $N^\delta \le H \le N^{1-\delta}$.
\endabstract

\endtopmatter

\document
\head
0. Introduction
\endhead
\adjustfootnotemark{-1}%
\footnotetext{Research supported by NSF grants DMS--0070720 and DMS--0244660.}%
\adjustfootnotemark{1}%
\footnotetext{Research supported by the American Institute of Mathematics and
NSF grants.}%

\noindent
Cram{\' e}r \cite{4} modeled the distribution of prime numbers by independent 
random variables $X_n$ (for $n\ge 3$) that take the value $1$ ($n$ is ``prime'') 
with probability $1/\log n$ and take the value $0$ ($n$ is ``composite'')  
with probability $1-1/\log n$.  If $p_n$ denotes the $n^{\text{th}}$ prime number 
this model predicts that 
$$
\lim_{N\to\infty}\frac1N\operatorname{card}
\{n:1\le n\le N, p_{n+1}-p_n > c\log p_n\} = e^{-c}
$$   
for all fixed positive real numbers $c$. Gallagher \cite{6} showed that
the above follows from Hardy \& Littlewood's \cite{10{\rm, p.\ 61}} quantitative version
of the prime $k$-tuple conjecture: If 
$\eusm D = \{d_1, d_2, \ldots, d_k\}$ is a set of $k$ distinct integers, then 
$$
\sum_{n\le x}\ \prod_{i=1}^k \Lambda(n+d_i)\ = \ (\eufm S(\eusm D)+o(1))x
\tag 1
$$
as $x\to\infty$ where $\eufm S(\eusm D)$ is the singular series
$$
\align
\eufm S(\eusm D)\ &= \sum\Sb q_1, \ldots, q_k\\ 1\le q_i < \infty\endSb
\bigg(\prod_{i=1}^k \frac{\mu(q_i)}{\phi(q_i)}\bigg)
\sum\Sb a_1, \ldots, a_k \\ 1\le a_i\le q_i \\ (a_i, q_i) = 1 \\ 
\sum a_i/q_i \in{\Bbb Z}\endSb  
e\Big(\sum_{i=1}^k \frac{a_i d_i}{q_i}\Big) \tag 2 \\
\intertext{where $e(\theta) = e^{2\pi i\theta}$. 
Hardy \& Littlewood showed that the right hand side above may be written more 
transparently as}
&= \ \prod_p \Big(1-\frac 1p\Big)^{\!-k} \Big(1-\frac{\nu_p(\eusm D)}{p}\Big)
\tag 3 \\
\endalign
$$
where $\nu_p(\eusm D)$ denotes the number of distinct residue classes 
modulo $p$ found among the members of $\eusm D$.  (See the remarks following
the proof of Lemma~3 in \S2.)  Since $\nu_p(\eusm D) = k$ for all sufficiently
large $p$, the product (3) is absolutely convergent.  Hence its value is
$0$ if and only if there is a prime $p$ for which $\nu_p(\eusm D) = p$.  
Gallagher \cite{6} showed that from (1) it follows that  
$$
\int_2^X (\psi(x+\lambda \log x)-\psi(x))^k\,dx \sim m_k(\lambda)X(\log X)^k
$$
when $\lambda\asymp 1$.  Here $m_k(\lambda) = {\Bbb E}(Y^k)$ is the 
$k^{\text{th}}$ moment of a Poisson random variable $Y$ with parameter 
$\lambda$, and $a \asymp b$ means that $a/b$ lies between two positive absolute constants.
Thus the distribution of $\pi(x+h)- \pi(h)$ is approximately
Poisson when $h\asymp \log X$, as predicted by the Cram{\' e}r model.  

%Colloquially we may express (1) as saying that the ``probability'' 
%of $n+d_1$, $\ldots$, $n+d_k$ all being prime is $\eufm S(\eusm D)/(\log n)^k$, 
%while the Cram{\' e}r model would predict the (obviously incorrect) 
%probability $1/(\log n)^{k}$.  The crux of Gallagher's argument is to 
%show that 
%$$
%\sum\Sb d_1,\ldots, d_k \\ 1\le d_i \le h \\ d_i \text{distinct}\endSb 
%\eufm S(\eusm D)\ \sim \ h^k
%\tag 4
%$$   
%as $h\to\infty$.  Thus, on average, $\eufm S(\eusm D) \sim 1$ and therefore 
%the distribution of primes in intervals of length $\asymp \log x$ 
%conforms to the predictions of the Cram{\' e}r model.  

In this paper we investigate the distribution of primes 
in longer intervals.  Let $H=H(N)$ be a function of $N$ such 
that $H=o(N)$ and $H/\log N \to \infty$ as $N\to \infty$. 
The Cram{\' e}r model predicts that the distribution of 
$\psi(n+H)-\psi(n)$ (for $n\le N$) is approximately normal with 
mean $\sim H$ and variance $\sim H\log N$.  Assuming a 
strong form of the Hardy-Littlewood conjecture (1) we will show that 
this prediction holds in the range where $H/\log N \to \infty$   
and $\log H/\log N \to 0$ as $N\to \infty$.  In the range $N^{\delta} 
\le H \le N^{1-\delta}$ 
we provide evidence showing that the Cram{\' e}r model is incorrect, and conjecture 
instead that the distribution of $\psi(n+H)-\psi(n)$ is approximately normal 
with mean $\sim H$ and variance $\sim H\log (N/H)$.   

When $h \asymp \log x$, the moments of $\psi(x+h)-\psi(x)$
and of $\psi(x+h)-\psi(x)-h$
are of the same order of magnitude.  However, for larger $h$ one would
expect that the moments of $\psi(x+h)-\psi(x)$ to be far larger than those of
$\psi(x+h)-\psi(x)-h$.  We obtain our conclusions on the distribution of $\psi(x+h) 
-\psi(x)$ by analyzing these latter, more delicate moments.  To facilitate this 
study, we 
set $\Lambda_0(n) = \Lambda(n) - 1$, with the result that
$$
\psi(x+h)- \psi(x) - h = \sum_{x<n\le x+h}\Lambda_0(n).
$$
Thus the main term is eliminated at the outset, which simplifies our calculations
considerably.  We recast (1) in
an equivalent form that pertains to $\Lambda_0(n)$:  If $d_1,\ldots, d_k$ are
distinct integers, then
$$
\sum_{n\le x}\ \prod_{i=1}^k \Lambda_0(n+d_i)\ = \ (\eufm S_0(\eusm D)+o(1))x
\tag 4
$$
as $x\to\infty$ where $\eufm S_0(\eusm D)$ is related to $\eufm S(\eusm D)$
by the identities
$$
\align
\eufm S_0(\eusm D) &= \sum_{\eusm I\subseteq \eusm D}
(-1)^{\text{card}\,\eusm I}\eufm S(\eusm I), \tag 5 \\
\eufm S(\eusm D) &= \sum_{\eusm I\subseteq \eusm D}\eufm S_0(\eusm I).
\tag 6
\endalign
$$
Here it is to be understood that $\eufm S_0(\emptyset)= \eufm S(\emptyset)=1$.  
One of the main steps in Gallagher's argument is to 
show that
$$
\sum\Sb d_1,\ldots, d_k \\ 1\le d_i \le h \\ d_i \text{distinct}\endSb 
\eufm S(\eusm D)\ \sim \ h^k
\tag 7
$$   
as $h\to\infty$.  There are $\big({h\atop k}\big)$ subsets $\eusm D$ under 
consideration, but each one occurs $k!$ times in the above sum.  Thus the above 
asserts that the mean value of $\eufm S(\eusm D)$ tends to $1$ as $h \to \infty$.  
Correspondingly, we need to estimate the quantities
$$
R_k(h) = \sum\Sb d_1,\ldots, d_k \\ 1 \le d_i \le h \\ d_i \text{distinct}\endSb 
\eufm S_0(\eusm D)\,.
\tag 8
$$
From (2) and (5) we see that
$$
\eufm S_0(\eusm D) = \sum\Sb q_1, \ldots, q_k\\ 1 < q_i < \infty\endSb
\bigg(\prod_{i=1}^k \frac{\mu(q_i)}{\phi(q_i)}\bigg)
\sum\Sb a_1, \ldots, a_k \\ 1\le a_i\le q_i \\ (a_i, q_i) = 1 \\ 
\sum a_i/q_i \in{\Bbb Z}\endSb  
e\Big(\sum_{i=1}^k \frac{a_i d_i}{q_i}\Big).
\tag 9
$$
The task of estimating averages of this expression is quite challenging,
but our burden is substantially lightened by work of Montgomery \& Vaughan \cite{16} 
concerning a strikingly similar quantity.  Let
$$
m_k(q;h) = \sum_{n=1}^q \Bigg(\sum \Sb m=1\\ (m+n,q)=1\endSb^h 
1\ - \ h\phi(q)/q\Bigg)^{\!\!k}
\tag 10
$$
be the $k^{\text{th}}$ centered moment of the number of reduced residues (mod $q$)
in an interval.  Lemma~2 of Montgomery \& Vaughan asserts that
$$
m_k(q;h) = q\Big(\frac{\phi(q)}q\Big)^{\!k}V_k(q;h)
\tag 11
$$
where
$$
V_k(q;h)= \sum\Sb d_1, \ldots, d_k \\ 1\le d_i \le h \endSb 
\sum\Sb q_1, \ldots, q_k\\ 1 < q_i|q\endSb
\bigg(\prod_{i=1}^k \frac{\mu(q_i)}{\phi(q_i)}\bigg)
\sum\Sb a_1, \ldots, a_k \\ 1\le a_i\le q_i \\ (a_i, q_i) = 1 \\
\sum a_i/q_i \in \Bbb Z \endSb  
e\Big(\sum_{i=1}^k \frac{a_i d_i}{q_i}\Big).
\tag 12
$$
When $k=1$, the conditions in the innermost sum cannot be fulfilled, and thus
$V_1(q;h)= 0$.  When $k=2$, the conditions in the innermost sum require that
$q_1 = q_2 = a_1 + a_2$.  Thus
$$
V_2(q;h) = \sum\Sb d|q \\ d > 1 \endSb\frac{\mu(d)^2}{\phi(d)^2}
\sum\Sb a=1 \\ (a,d)=1 \endSb^d |E(a/d)|^2
\tag 13
$$ 
where 
$$
E(\alpha) = \sum_{m=1}^h e(m\alpha).
\tag 14
$$
Montgomery \& Vaughan showed that
$$
V_k(q;h) \ \ll_k \ (hq/\phi(q))^{k/2}
\big(1 + O\big(h^{-1/(7k)}(q/\phi(q))^{2^k+k/2}\big)\big)
$$
for each positive integer $k$.  Unfortunately, this is not quite sharp enough
for our present purposes, so our first job is to refine the above.

\proclaim{Theorem 1}In the above notation,
$$
V_k(q;h) = \mu_k V_2(q;h)^{k/2}
+O_k\bigg(h^{k/2-1/(7k)}\Big(\frac q{\phi(q)}\Big)^{\!2^k+k/2}\bigg)
\tag 15
$$
for every positive integer $k$, where $\mu_k= 1\cdot3\cdots(k-1)$
if\/ $k$ is even, and $\mu_k = 0$ if\/ $k$ is odd.
\endproclaim

    Here the main term is the $k^{\text{th}}$ moment of a normal random variable
with expectation $0$ and variance $V_2(q;h)$.  We remark that the work of 
Granville \& Soundararajan [7] (see \S 6a) places restrictions on the 
uniformity (in $k$) with which (15) can possibly hold. 
With Theorem 1 in hand, 
we are able to estimate the $R_k(h)$.
        
\proclaim{Theorem 2}Let $h$ be an integer, $h > 1$, and suppose that $R_k(h)$ 
is defined as in {\rm(8)}.  Then
$$
R_k(h) = \mu_k(-h\log h + Ah)^{k/2} +O_k\big( h^{k/2-1/(7k)+\varepsilon}\big)
$$
for any nonnegative integer $k$, where $A = 2 - C_0 - \log 2\pi$ and $C_0$ denotes
Euler's constant.
\endproclaim

    For the smallest values of $k$, one can be more precise, since it is clear
that $R_0(h)=1$, and that $R_1(h) = 0$.  Also, from (5) and (48) we know that
$$
R_2(h) = -h\log h + Ah +O(h^{1/2+\varepsilon}).
\tag 16
$$      
From (6) it follows that the left hand side of (7) is
$$
\sum_{r=0}^k\Big({k\atop r}\Big)R_r(h)(h-r)(h-r-1)\cdots(h-k+1).
$$      
Hence we obtain Gallagher's estimate (7) in the more precise form
$$
\sum\Sb d_1,\ldots, d_k\\ 1 \le d_i \le h \\ d_i \text{distinct}\endSb 
\eufm S(\eusm D)
= h^k -\Big({k\atop 2}\Big)h^{k-1}\log h 
+ \Big({k\atop 2}\Big)(1-C_0-\log 2\pi)h^{k-1} +O\big(h^{k-3/2+\varepsilon}\big).
\tag 17
$$

     We put
$$
\align
M_K(N;H) &= \sum_{n=1}^N (\psi(n+H)-\psi(n)-H)^K,
\tag 18 \\
\intertext{and note that this is}
&= \sum_{n=1}^N\Big(\sum_{h=1}^H \Lambda_0(n+h)\Big)^{\!K}
= \sum\Sb h_1,\ldots, h_K\\ 1\le h_i\le H\endSb
\ \sum_{n=1}^N\ \prod_{i=1}^K \Lambda_0(n+h_i).
\tag 19
\endalign
$$
Here the $h_i$ are not necessarily distinct, but once the distinct values have
been identified, and their multiplicities accounted for, we can appeal to 
(4).  The quantities $R_k(h)$ arise in the main term.  Thus from Theorem 2
we can derive an asymptotic estimate for the above, provided that the error
term in (4) is sufficiently small and $H$ is not too large.

\proclaim{Theorem 3}Let $E_k(x;\eusm D)$ be defined by the relation
$$
\sum_{n\le x}\ \prod_{i=1}^k \Lambda(n+d_i)\ = \ \eufm S(\eusm D)x +E_k(x;\eusm D),
$$
and suppose that 
$$
E_k(x;\eusm D)\ll N^{1/2+\varepsilon}
\tag 20
$$ 
uniformly for $1\le k\le K$, $0\le x\le N$, and distinct $d_i$ satisfying 
$1\le d_i\le H$.  Then
$$
\align
M_K(N;H) &= \mu_KH^{K/2}\int_1^N(\log x/H + B)^{K/2}\,dx 
\\
&+O\Big(N(\log N)^{K/2} H^{K/2} \Big(\frac{H}{\log N}\Big)^{-1/(8K)} + 
H^{K} N^{1/2+\varepsilon} \Big)
\tag 21
\\
\endalign
$$
uniformly for $\log N \le H\le N^{1\!/K}$, where 
$B = 1 - C_0-\log 2\pi$
and $C_0$ denotes Euler's constant.
\endproclaim

     In the case $k=1$, the set $\eusm D$ is a singleton, $\eufm S(\eusm D) = 1$,
and the hypothesis that $E_1(N;\{1\})\ll N^{1/2+\varepsilon}$ is equivalent to
the Riemann Hypothesis (RH). 

In place of (20) if we assume only that $E_k(x;\eusm D) \ll E$ for 
some $E \ge N^{1/2+\varepsilon}$ then (21) holds with the modified error term 
$$
 O\Big(N(\log N)^{K/2} H^{K/2} \Big(\frac{H}{\log N}\Big)^{\!-1/(8K)} + 
H^{K} E \Big).
$$       

     We note that in deriving Theorem 3 from Theorem 2, we start with $h_i$ that 
are not necessarily distinct, and must reduce to distinct $d_i$, which is the 
reverse of the problem encountered in deriving Theorem 2 from Theorem 1, where 
we start with distinct $d_i$, and want to appeal to an estimate involving not 
necessarily distinct $m_i$.      

Since the $\mu_K$ are the moments of a normal random variable 
with mean $0$ and variance $1$ we deduce from Theorem 3 the 
following Corollary. 

\proclaim{Corollary 1} Let $H=H(N)$ be a function of $N$ such that 
$$
\frac H{\log N} \to +\infty, \qquad \frac{\log H}{\log N} \to 0
$$
as $N\to \infty$.  
Assume that the hypothesis {\rm(20)} holds for arbitrarily large $K$.  
Then the distribution of $\psi(n+H)-\psi(n)-H$ for $n\le N$ is 
approximately normal with mean $0$ and variance $H \log N$, in the sense that
the proportion of $n\le N$ for which $\psi(n+H)-\psi(n)-H\le cH\log N$
tends to $\Phi(c)$ as $N \to \infty$, uniformly for $|c|\le C$.  Here 
$\Phi(u) =\frac1{\sqrt{2\pi}}\int_{-\infty}^u e^{-v^2\!/2}\,dv$ is the cumulative
distribution function of a normal random variable with mean $0$ and variance $1$.
\endproclaim

     For larger $H$, Theorem 3 furnishes only a limited 
number of moments and we cannot deduce a distribution result.  
However we expect that the contributions of the $E_k(x;\eusm D)$ cancel sufficiently 
so as not to overwhelm the main term: 

\proclaim{Conjecture 1}For each positive integer $K$,
$$
M_K(N;H) = (\mu_K+o(1))N\Big(H\log \frac NH\Big)^{\!K/2}
$$
uniformly for $(\log N)^{1+\delta}\le H \le N^{1-\delta}$.

\endproclaim

    This implies the weaker 
         
\proclaim{Conjecture 2}Suppose that $(\log N)^{1+\delta}\le H \le N^{1-\delta}$.  
The distribution of $\psi(x+H)-\psi(x)$ for $0\le x\le N$ is approximately normal 
with mean $H$ and variance $H\log N/H$.
\endproclaim

     Certainly Conjecture 2 does not hold when $H \asymp N$, but 
perhaps it holds whenever $H=o(N)$.  It would be interesting to investigate 
more thoroughly what happens in this range.  
   
     Hardy \& Littlewood \cite{10} provided heuristics that point toward the
quantitative prime $k$-tuple conjecture (1).  In \S4 we argue in the same spirit
to obtain indications in favor of Conjecture~1. 
         
   To obtain further support for our conjectures, we interpret the situation
in terms of the zeros of the Riemann zeta function.  We recall that Goldston \&
Montgomery \cite{8} showed that if RH is true, then the stronger form 
($F(\alpha) \sim 1$) of the Pair Correlation Conjecture as formulated by 
Montgomery \cite{13} is equivalent to the case $K=2$ of the Conjecture above.  
In the same spirit, Chan \cite{3} has shown (assuming RH) that Conjecture 1 is 
equivalent to the assertion that
$$
\int_1^X\Big(\sum_{0<\gamma\le T}\cos( \gamma \log x)\Big)^{\!k}\,dx 
=(\mu_k +o(1))X\Big(\frac T{4\pi}\log T\Big)^{k/2}\,.
\tag 22
$$
Viewed in this way, we see that the Pair Correlation Conjecture asserts that the 
variance of the sum
$$
\sum_{0<\gamma\le T}\cos(\gamma\log x)
$$
is the same as it would be if it were a sum of uncorrelated random variables, and 
Conjecture~1 asserts that this same sum has the same normal distribution that it would
have if the terms were independent random variables.  In somewhat the same vein,
Bogomolny \& Keating \cite{1} used Hardy--Littlewood conjectures concerning  primes
to arrive at the $n$ level correlation function of zeros of the zeta function.  

     Freeman Dyson observed that the Pair Correlation Conjecture is analogous to 
known properties of the spacings of the eigenvalues of certain families of random matrices.  
We note that (22) has a similar analogue in random matrix theory.  Let $U(N)$ denote the
classical compact group of unitary $N\times N$ matrices.  For $A\in U(N)$, let
$e(\theta_1), \ldots, e(\theta_N)$ denote the eigenvalues of $A$.  Rains \cite{19}
has observed that if $M$ is an integer, $|M|\ge N$, then the point 
$(M\theta_1, \ldots, M\theta_N)$ is exactly uniformly distributed in $\Bbb T^N$
as $A$ varies with respect to the Haar measure $dA$ on $U(N)$.  It follows in 
particular that the distribution of
$$
\operatorname{Re\ Trace} A^M = \sum_{n=1}^N \cos(2\pi M\theta_n)
$$
is exactly the same as the distribution of
$$
\sum_{n=1}^N \cos 2\pi X_n
$$
where where the $X_n$ are independent random variables, each one uniformly distributed
on $[0,1]$.  It follows by the Central Limit Theorem that this distribution tends
toward a normal distribution with mean $0$ and variance $N/2$.  By an easy calculation
it can also be shown that the $k^{\text{th}}$ moment of this distribution is $0$ if $k$ 
is odd, and is $\sim \mu_k(N/2)^{k/2}$ if $k$ is even.
 
     As for numerical studies, Brent~\cite{2} has compiled evidence not only for
(1) but also for the stronger hypothesis (20).  Odlyzko~\cite{18}
and Forrester \& Odlyzko~\cite{5} have found that the local distribution of
the zeros of the zeta function fits well with predictions based on random matrix
theory.  The authors~\cite{15} have reported on numerical evidence in support of
the conjectures.  Finally, Chan~\cite{3{\rm; pp.\ 36, 49, 63}} has assembled evidence 
in favor of~(22).        

%     Let $p_n$ denote the $n^{\text{th}}$ prime number.
%Cram\'er modeled the primes by constructing a probability space in which 
%an integer $n\ge 3$  is `prime' with probability $1/\log n$.  This model 
%suggests that 
%$$
%\lim_{N\to\infty}\frac1N\operatorname{card}
%\{n:1\le n\le N, p_{n+1}-p_n > c\log p_n\} = e^{-c}
%$$   
%for all fixed positive real numbers $c$, and Gallagher \cite{7} showed that
%the above follows from (1).  Concerning somewhat
Cram\'er's model suggests that
$$
\pi(x+(\log x)^a) - \pi(x)\sim (\log x)^{a-1}
\tag 23
$$
as $x\to \infty$ with $a$ fixed, $a > 2$.  This, however, is known to be false,
since Maier \cite{11} showed that 
$$
\mathop{\overline{\underline{\lim}}}_{x\to\infty}
\frac{\pi(x+(\log x)^a)-\pi(x)}
{(\log x)^{a-1}}\  \gtrless \ 1
$$
for any fixed $a > 0$ (for general results of 
this nature see Granville \& Soundararajan [7]).  
Presumably (23) is valid for most $x$, and the exceptions
discovered by Maier are quite rare. Indeed Selberg [18] showed that on RH, (23) 
holds if $a>2$ for almost all $x$, and Corollary 1 shows on hypothesis (20) 
that (23) holds if $a>1$ for almost all $x$.  As for longer intervals, suppose that
$\alpha$ is fixed, $0 < \alpha < 1$.  Cram\'er's 
model would predict that $\psi(x+x^\alpha)-\psi(x)- x^\alpha$ is 
approximately normally distributed with mean $0$ and variance $X^\alpha\log X$
as $x$ runs over the range $X\le x\le 2X$.  Our Conjecture 1 predicts normal 
distribution, but with a variance that is smaller by a factor of $1-\alpha$.  Thus
it seems that in this range, Cram\'er's model is not just occasionally inaccurate,
but instead is actually inaccurate on average.

\head 
1. Proof of Theorem 1
\endhead

\noindent
Montgomery \& Vaughan \cite{16} devised a useful basic inequality (their Lemma~1),
which we now quote.

\proclaim{Lemma 1}Let $r_1,\ldots, r_k$ be squarefree integers, set 
$r = [r_1, \ldots, r_k]$, and suppose that any prime dividing $r$ divides at least
two of the $r_i$.  Then for any complex-valued functions $G_1, \ldots, G_k$ defined
on $(0,1]$ we have
$$
\bigg|\sum\Sb b_1, \ldots, b_k \\ 1\le b_i \le r_i \\ 
\sum b_i/r_i \in\Bbb Z \endSb \prod_{i=1}^k G_i(b_i/r_i)\bigg|
\le \frac 1r\prod_{i=1}^k\Big(r_i\sum_{b_i=1}^{r_i}|G_i(b_i/r_i)|^2\Big)^{\!1/2}.
$$  
\endproclaim

     Montgomery \& Vaughan \cite{17} have derived several variants of the above; an
exposition of such variants is found in Chapter~8 of Montgomery \cite{14}.  For our
present purposes a different type of variant is useful.

\proclaim{Lemma 2}Let $q_1, \ldots, q_k$ be squarefree integers, each one strictly
greater than $1$, and put $d = [q_1, \ldots, q_k]$.  Let $G$ be a complex-valued 
function defined on $(0,1)$, and suppose that $G_0$ is a nondecreasing function 
on the positive integers such that
$$
\sum_{a=1}^{q-1} |G(a/q)|^2 \le qG_0(q)
\tag 24
$$
for all squarefree integers $q>1$.  Then
$$
\bigg|\sum\Sb a_1, \ldots, a_k \\ 0 < a_i < q_i \\ 
\sum a_i/q_i \in\Bbb Z \endSb \prod_{i=1}^k G(a_i/q_i)\bigg|
\le \frac 1d\prod_{i=1}^k q_iG_0(q_i)^{1/2}.
$$
\endproclaim

\demo{Proof}We write $q_i= r_is_i$ where the $s_i$ are pairwise relatively prime
and any prime dividing $[r_1, \ldots, r_k]$ divides at least two of the $r_i$.
That is, $r_i = (q_i, \prod_{j\ne i}q_j)$.  Clearly $d = rs_1\cdots s_k$ where
$r = [r_1,\ldots, r_k]$.
The condition $\sum a_i/q_i\in \Bbb Z$ forces $s_i|a_i$ for all $i$.  Hence, on
writing $a_i = s_ib_i$, we find that the left hand side above is
$$
= \bigg|\sum\Sb b_1, \ldots, b_k \\ 0 < b_i < r_i \\ \sum b_i/r_i\in\Bbb Z\endSb
\prod_{i=1}^k G(b_i/r_i)\bigg|\,.
$$
If there is an $i$ for which $r_i=1$, then the conditions in the above sum cannot
be fulfilled, the sum is empty, and there is nothing to prove.  Thus we may assume
that $r_i > 1$ for all~$i$.
In Lemma 1 we take $G_i(x)=G(x)$ for $0 < x < 1$, and $G(1)=0$.  Thus by Lemma~1
and the hypothesis (24), the above is
$$
\le \frac1r\prod_{i=1}^k\Big(r_i\sum_{b_i=1}^{r_i-1}r_i|G(b_i/r_i)|^2\Big)^{\!1/2}
\le \frac1r\prod_{i=1}^k \big(r_i^2G_0(r_i)\big)^{\!1/2}.
$$
Since $G$ is nondecreasing, the above is
$$
\le \frac1r\prod_{i=1}^k \Big(r_i^2G_0(q_i)\big)^{\!1/2}
=\frac1d\prod_{i=1}^k \Big(q_i^2G_0(q_i)\big)^{\!1/2},
$$
as desired.
\enddemo

    We now begin the main body of the proof of Theorem 1.  We take $k$ to be fixed,
so that the dependence of implicit constants on $k$ is suppressed.  If $k$ is odd,
then the desired estimate is already found in (18) of 
Montgomery \& Vaughan \cite{16}.  Thus we may assume that $k$ is even.  
From (12) and (14) it is clear that
$$
V_k(q;h) = \sum\Sb q_1, \ldots, q_k\\ 1 < q_i|q\endSb
\bigg(\prod_{i=1}^k \frac{\mu(q_i)}{\phi(q_i)}\bigg)
\sum\Sb a_1, \ldots, a_k \\ 1\le a_i\le q_i \\ (a_i, q_i) = 1 \\
\sum a_i/q_i \in \Bbb Z \endSb  
\prod_{i=1}^k E(a_i/q_i).
$$
In Lemmas 7 and 8 of Montgomery \& Vaughan \cite{16}, it is shown that all 
contributions to the above are
$$
\ll h^{k/2-1/(7k)}\Big(\frac q{\phi(q)}\Big)^{\!2^k+k/2},
$$
except for those terms for which the $q_i$ are equal in pairs, with no further 
equalities among the $q_i$.  There are $(k-1)(k-3)\cdots3\cdot1 = \mu_k$ ways in 
which this pairing can occur.  Take the pairing to be $q_i = q_{k/2+i}$, and set 
$b_i = a_i+a_{k/2+i}$.  Thus the terms that remain to be estimated are precisely
$$
\mu_k\sum\Sb q_1,\ldots, q_{k/2} \\ 1 < q_i|q \\ q_i \text{ distinct}\endSb\
\prod_{i=1}^{k/2}\frac{\mu(q_i)^2}{\phi(q_i)^2}
\sum\Sb b_1, \ldots, b_{k/2} \\ 1\le b_i \le q_i \\ \sum b_i/q_i \in\Bbb Z \endSb
\prod_{i=1}^{k/2} J(b_i, q_i)
\tag 25
$$
where 
$$
J(b,r) = \sum\Sb a=1 \\ (a,r)=1 \\ (b-a,r)=1\endSb^r E\big(\frac ar\Big)
E\Big(\frac{b-a}r\Big).
\tag 26
$$
First we show that the condition that the $q_i$ should be distinct in (25) can be
dropped.  To see this, put $F(\alpha) = \min (h, 1/\|\alpha\|)$ where 
$\|\theta\| = \min_{n\in\Bbb Z}|\theta-n|$ is the distance from $\theta$ to the
nearest integer.  Thus 
$$
|E(\alpha)|\le F(\alpha)
\tag 27
$$
for all $\alpha$.  Let
$\eusm Q$ denote the set of those $k$-tuples $(q_1, \ldots, q_k)$ such that
$1 < q_i|q$ for all $i$, and with the property that among the $q_i$ there are
three or more of them that are equal.  In proving their Lemma~8 (see the treatment
of $T_3$), Montgomery \& Vaughan \cite{16} establish that
$$
\sum_{\boldkey q\in\eusm Q}
\bigg(\prod_{i=1}^k \frac{|\mu(q_i)|}{\phi(q_i)}\bigg)
\sum\Sb a_1, \ldots, a_k \\ 1\le a_i\le q_i \\ (a_i, q_i) = 1 \\
\sum a_i/q_i \in \Bbb Z \endSb  
\prod_{i=1}^k F(a_i/q_i)\
\ll \ h^{k/2-1/(7k)}\Big(\frac q{\phi(q)}\Big)^{\!2^k+k/2}.
$$
Since this majorizes the difference between (25) and
$$
\mu_k\sum\Sb q_1,\ldots, q_{k/2} \\ 1 < q_i|q \endSb\
\prod_{i=1}^{k/2}\frac{\mu(q_i)^2}{\phi(q_i)^2}
\sum\Sb b_1, \ldots, b_{k/2} \\ 1\le b_i \le q_i \\ \sum b_i/q_i \in\Bbb Z \endSb
\prod_{i=1}^{k/2} J(b_i, q_i),
\tag 28
$$
it follows that we can continue with the above expression.  Suppose that 
$0 < b_i < q_i$ for exactly $j$ values of $i$, and that $b_i = q_i$ for the 
remaining $k/2-j$ values of $i$.  Since there are $\big({k/2\atop j}\big)$ ways 
of choosing the $j$ indicies, we see that the above is
$$
\mu_k\sum_{j=0}^{k/2}\Big({k/2\atop j}\Big)V_2(q;h)^{k/2-j}W_j(q;h)
\tag 29
$$
where $W_0(q;h)=1$ and
$$
W_j(q;h) = \sum\Sb q_1, \ldots, q_j \\ 1 < q_i|q\endSb\
\prod_{i=1}^j\frac{\mu(q_i)^2}{\phi(q_i)^2}
\sum\Sb b_1, \ldots, b_j \\ 0 < b_i < q_i \\ \sum b_i/q_i \in\Bbb Z\endSb
\prod_{i=1}^j J(b_i,q_i).
\tag 30
$$
Here the term $j=0$ gives the desired main term.  Thus it remains to show that the
other terms are smaller.  

     To prepare for an application of Lemma~2, we estimate $J(b,r)$.  By (27) we see 
that if $0 < b \le r/2$ and $r < h$, then
$$
\align
J(b,r) &\ll  \sum_{b/2<a<b}\frac{r^2}{b(b-a)}\
+\sum_{b < a \le 3b/2}\frac{r^2}{b(a-b)}\ +  \sum_{3b/2<a\le 2r/3}\frac{r^2}{a^2} \\
&\ll \frac{r^2}b\log 2b.
\tag 31
\endalign
$$
Here half the ranges of $a$ have been omitted, since by symmetry they contribute the
same amount as the listed sums.
Similarly, if $0 < b \le r/h$ and $r\ge h$, then
$$
J(b,r) \ll \sum_{0<a\le 2r/h}h^2 \ + \sum_{2r/h<a\le r/2}\frac{r^2}{a^2} \ll rh\,. 
\tag 32
$$
Finally, if $r/h < b\le r/2$ and $r\ge h$, then
$$
\align
J(b,r) &\ll \sum_{b/2<a\le b-r/h}\frac{r^2}{b(b-a)}\ 
+ \sum_{b-r/h< a < b}\frac{rh}a\
+ \sum_{b < a \le b + r/h}\frac{rh}a \\ 
&\qquad + \sum_{b+r/h<a\le 4b/3}\frac{r^2}{b(a-b)}\
+ \sum_{4b/3<a\le 2r/3}\frac{r^2}{a^2} \\
&\ll \frac{r^2}b\log (2bh/r).
\tag 33
\endalign
$$
From (31) we see that if $r < h$, then
$$
\sum_{0<b<r}J(b,r)^2 \ll r^4,
$$
and from (32) and (33) we see that if $r\ge h$, then
$$
\sum_{0<b<r}J(b,r)^2 \ll r^3h.
$$
Altogether,
$$
\sum_{0 < b < r}J(b,r)^2 \ll r^3\min(r, h).
\tag 34
$$
On taking $G_0(r)= Chr^2$ in Lemma 2, we find that
$$
\sum\Sb b_1, \ldots, b_j \\ 0 < b_i < q_i \\ \sum b_i/q_i \in \Bbb Z\endSb
\prod_{i=1}^j J(b_i,q_i) \ll \frac1d\prod_{i=1}^j \big(q_i^2h^{1/2}\big),
$$
and hence
$$
W_j(q;h) \ll h^{j/2}\sum_{d|q}\frac1d
\Big(\sum_{r|d}\frac{\mu(r)^2r^2}{\phi(r)^2}\Big)^{\!j}
= h^{j/2}\prod_{p|q}\Big(1+\frac1p\Big(1+\frac{p^2}{(p-1)^2}\Big)^{\!j}\Big)
\ll h^{j/2}\Big(\frac q{\phi(q)}\Big)^{\!2^j}.
$$
To apply this in (29), we need also a bound for $V_2(q;h)$.  To this end we 
note that
$$
\align
V_2(q;h) &\le \sum_{d|q}\frac{\mu(d)^2}{\phi(d)^2}\sum_{a=1}^{d-1} F(a/d)^2
\ll h\sum_{d|q}\frac{\mu(d)^2d}{\phi(d)^2} \\
&= h\prod_{p|q}\Big(1+\frac p{(p-1)^2}\Big) \ll h\frac q{\phi(q)}\,.
\endalign
$$
(By a different method it can be shown that $V_2(q;h)\le hq/\phi(q)$.  See Hausman
\& Shapiro \cite{11} and (3) of Montgomery \& Vaughan \cite{16}.)  From (30) we
see that $W_1(q;h)=0$, since the inner sum is empty.  On applying the above estimates
for $2\le j\le k/2$, we see that the expression (29) is
$$
\mu_kV_2(q;h)^{k/2} +O\big(h^{k/2-1}(q/\phi(q))^{2^{k/2}}\big).
$$
Here the error term is majorized by that in (15), so the proof is complete.

\head 
2. Proof of Theorem 2
\endhead
\noindent
We begin with two lemmas.

\proclaim{Lemma 3}{\rm(Hardy--Littlewood)} Let
$$
A(q_1, \ldots, q_k) = \sum\Sb a_1, \ldots, a_k \\ 1\le a_i\le q_i \\
(a_i,q_i)=1 \\ \sum a_i/q_i \in\Bbb Z \endSb 
e\Big(\sum_{i=1}^k \frac{d_ia_i}{q_i}\Big). 
$$
If\/ $q_i = q_i'q_i''$ with $(\prod q_i',\prod q_i'')=1$, then
$$
A(q_1, \ldots, q_k) = A(q_1', \ldots, q_k')A(q_1'', \ldots, q_k'').
\tag 35
$$
For any prime number $p$,
$$
\sum\Sb q_1, \ldots, q_k \\ q_i|p \endSb\ \prod_{i=1}^k \frac{\mu(q_i)}{\phi(q_i)}
A(q_1, \ldots, q_k) = \Big(1-\frac1p\Big)^{\!-k}\Big(1-\frac{\nu_p(\eusm D)}p\Big)
\tag 36
$$
where $\nu_p(\eusm D)$ is the number of distinct residue classes modulo $p$
found among the members of\/ $\eusm D = \{d_1, \ldots, d_k\}$.  Finally,
$$
\sum\Sb q_1, \ldots, q_k \\ 1\le q_i<\infty\endSb\
\prod_{i=1}^k \frac{\mu(q_i)^2}{\phi(q_i)}|A(q_1,\ldots,q_k)| < \infty.
\tag 37
$$
\endproclaim

\demo{Proof}We follow the argument of Hardy \& Littlewood 
\cite{10{\rm, pp. 56--61}}, but with
some helpful amplifications.  We write
$$
\frac{a_i}{q_i} \equiv \frac{a_i'}{q_i'}+\frac{a_i''}{q_i''}\pmod 1.
$$
By the Chinese Remainder Theorem, each reduced residue $a_i$ modulo $q_i$ corresponds
to a pair $a_i',a_i''$ of reduced residues modulo $q_i'$ and $q_i''$, respectively.
Also, $\sum a_i/q_i \in\Bbb Z$ if and only if $\sum a_i'/q_i'\in\Bbb Z$ and
$\sum a_i''/q_i''\in\Bbb Z$.  This gives (35).

    If each $q_i$ is either $1$ or $p$, then
$$
\sum_{r=1}^p e\Big(\sum_{i=1}^k \frac{a_ir}{q_i}\Big)
=\cases p &\text{if }\sum_{i=1}^k a_i/q_i \in \Bbb Z, \\
        0 &\text{otherwise.}
\endcases               
$$      
Thus the left hand side of (36) is 
$$
\align
&= \frac1p\sum_{r=1}^p\sum\Sb q_1,\ldots,q_k \\ q_i|p\endSb\
\prod_{i=1}^k\frac{\mu(q_i)}{\phi(q_i)}
\sum\Sb a_1,\ldots,a_k \\1\le a_i\le q_i\\ (a_i,q_i)=1\endSb
e\Big(\sum_{i=1}^k\frac{a_i(d_i-r)}{q_i}\Big) \\
&= \frac1p\sum_{r=1}^p\prod_{i=1}^k\Big(1 - \frac1{p-1}
\sum_{0<a<p}e\Big(\frac{a(d_i-r)}p\Big)\Big).
\endalign
$$
Here the innermost sum is $p-1$ or $-1$, according as $r\equiv d_i\pmod p$, or
not.  Thus if $r\equiv d_i\pmod p$, then this factor of the product is $0$.
There are $\nu_p(\eusm D)$ such values of $r$.  For the remaining $p-\nu_p(\eusm D)$
values of $r$, each factor of the product is $p/(p-1)$.  Hence the above is
$$
= \frac{p - \nu_p(\eusm D)}p\Big(\frac p{p-1}\Big)^{\!k},
$$
which gives (36).

    From (35) we see that
$$
\sum\Sb q_1,\ldots, q_k \\ q_i|Q\endSb \prod_{i=1}^k\frac{\mu(q_i)^2}{\phi(q_i)}
|A(q_1, \ldots, q_k)|
= \prod_{p|Q}\Big(\sum\Sb q_1,\ldots, q_k \\ q_i|p\endSb 
\prod_{i=1}^k\frac1{\phi(q_i)}
|A(q_1, \ldots, q_k)|\Big).
\tag 38
$$      
Put $D = \prod_{i<j}(d_j-d_i)$.  For primes $p|D$ we make no attempt to simplify
the above expression.  However, there are only finitely many such primes, and
for $p\nmid D$, the $d_i$ are distinct modulo $p$. For such primes we
evaluate the factor more explicitly.  Let $\eusm J\subseteq\{1,\ldots,k\}$ with
$j = \operatorname{card}\eusm J$, suppose that $q_i=p$ for $i\in\eusm J$,
$q_i = 1$ for $i\notin\eusm J$, and that the $d_i$ are distinct modulo $p$.  
Then
$$
A(q_1, \ldots, q_k) = \frac1p\sum_{r=1}^p
\sum\Sb a_1, \ldots, a_k \\ 1\le a_i\le q_i\\ (a_i,q_i)=1\endSb
e\Big(\frac{a_i(d_i-r)}{q_i}\Big) 
= \frac1p\sum_{r=1}^p\ \prod_{i=1}^k\Big(
\sum\Sb 1\le a_i\le q_i\\(a_i,q_i)=1\endSb e\Big(\frac{a_i(d_i-r)}{q_i}\Big)\Big).
$$
If $i\notin\eusm J$, then the innermost sum is $1$ for all $r$.  If
$i\in\eusm J$, then the innermost sum is $p-1$ if $r\equiv d_i\pmod p$, and 
$-1$ otherwise.  Thus there are $j$ values of $r$ for which one factor is $p-1$,
$j-1$ factors are $-1$, and all other factors are $1$.  For the remaining $p-j$ 
values of $r$, there are $j$ factors that are $-1$ and the remaining factors are 
$1$.  Thus the above is
$$
= \frac1p\big(j(p-1)(-1)^{j-1}+(p-j)(-1)^j\big) = (-1)^{j-1}(j-1). 
$$
Hence $|A(q_1, \ldots, q_k)| = |j-1|$, so it follows that the expression (38) is
$$
\prod\Sb p|Q\\ p|D\endSb\bigg(\sum\Sb q_1,\ldots, q_k \\ q_i|p\endSb 
\prod_{i=1}^k\frac1{\phi(q_i)}
|A(q_1, \ldots, q_k)|\bigg) 
\times
\prod\Sb p|Q\\ p\nmid D\endSb
\bigg(1+\sum_{j=2}^k\Big({k\atop j}\Big)\frac{j-1}{(p-1)^j}\bigg).
$$
Since this last product converges when extended over all primes, we have (37),
and the proof is complete.
\enddemo

     From (35) and (36) we see that
$$
\sum\Sb q_1, \ldots, q_k \\ q_i|Q \endSb \prod_{i=1}^k \frac{\mu(q_i)}{\phi(q_i)}
A(q_1, \ldots, q_k) 
= \prod_{p|Q}\Big(1-\frac1p\Big)^{\!-k}\Big(1-\frac{\nu_p(\eusm D)}p\Big)
\tag 39
$$
for any positive integer $Q$.  By (37) it follows that
$$
\align
\eufm S(\eusm D) &= \lim_{y\to\infty}
\sum\Sb q_1, \ldots, q_k \\ p|q_i\Rightarrow p\le y\endSb
\prod_{i=1}^k\frac{\mu(q_i)}{\phi(q_i)}A(q_1,\ldots, q_k),
\intertext{which by (39) is}
&= \lim_{y\to\infty}
\prod_{p\le y}\Big(1-\frac1p\Big)^{\!-k}\Big(1-\frac{\nu_p(\eusm D)}p\Big).
\endalign
$$
Thus the expressions (2) and (3) are equal.  

    Suppose that $1\le d_i\le h$ for all $i$.  Then $\nu_p(\eusm D)=k$ for
all primes $p>h$, and thus if $y\ge h$, then 
$$
\prod_{p>y}\Big(1-\frac1p\Big)^{\!-k}\Big(1-\frac{\nu_p(\eusm D)}p\Big)
= \prod_{p>y}\Big(1 +O_k\Big(\frac1{p^2}\Big)\Big) =1+O_k\Big(\frac1{y\log y}\Big).
\tag 40
$$
Since $\nu_p(\eusm D)\ge 1$ for all $p$, it is evident that
$$
\prod_{p\le h}\Big(1-\frac1p\Big)^{\!-k}\Big(1-\frac{\nu_p(\eusm D)}p\Big)
 \ll_k (\log h)^{k-1}.
\tag 41
$$
On combining this with (40), we see that
$$
\eufm S(\eusm D) \ll (\log h)^{k-1}.
\tag 42
$$
From (5) it follows additionally that 
$$
\eufm S_0(\eusm D) \ll (\log h)^{k-1}.
\tag 43
$$
From (39)--(41) we see that if $1\le d_i \le h$ for all $i$ and $y\ge h$, then
$$
\eufm S(\eusm D) = \sum \Sb q_1,\ldots, q_k \\ p|q_i\Rightarrow p\le y\endSb
\prod_{i=1}^k\frac{\mu(q_i)}{\phi(q_i)}A(q_1,\ldots, q_k)\ 
+ \ O_k\Big(\frac{(\log y)^{k-2}}y\Big). 
\tag 44
$$ 
By combining this with (5), we see also that
$$
\eufm S_0(\eusm D) 
= \sum \Sb q_1,\ldots, q_k \\ q_i > 1 \\ p|q_i\Rightarrow p\le y\endSb
\prod_{i=1}^k\frac{\mu(q_i)}{\phi(q_i)}A(q_1,\ldots, q_k)\ 
+ \ O_k\Big(\frac{(\log y)^{k-2}}y\Big). 
\tag 45
$$

\proclaim{Lemma 4}Let $E(\alpha)$ be defined as in {\rm(14)}.  If\/ $q$ is 
divisible by every prime number $p\le h^2$, then
$$
\sum_{d|q}\frac{\mu(d)^2}{\phi(d)^2}\Big(\sum\Sb a=1\\(a,d)=1\endSb^d
|E(a/d)|^2 \ - \phi(d)h\Big)\ = \ h^2 - h\log h + Bh+O\big(h^{1/2+\varepsilon}\big)
\tag 46
$$
where $B = 1- C_0-\log 2\pi$ and $C_0$ denotes Euler's constant.
\endproclaim

\demo{Proof}Since $|E(\alpha)|^2= \sum_{|m|\le h}(h-|m|)e(m\alpha)$, it follows that
$$
\sum\Sb a=1\\(a,d)=1\endSb^d |E(a/d)|^2 = \sum_{|m|\le h}(h-|m|)c_d(m)
$$
where $c_d(m)$ is Ramanujan's sum.  Now $c_d(0)=\phi(d)$, and $c_d(-m)=c_d(m)$,
so the above is 
$$
\phi(d)h + 2\sum_{m=1}^h (h-m)c_d(m).
$$
Hence the left hand side of (46) is
$$
= 2\sum_{m=1}^h(h-m)\sum_{d|q}\frac{\mu(d)^2}{\phi(d)^2}c_d(m).
$$
Here the sum over $d$ is
$$
\prod_{p|q}\Big(1+\frac{c_p(m)}{(p-1)^2}\Big)
= \prod\Sb p|q\\ p|m\endSb \Big(1+\frac1{p-1}\Big)\prod\Sb p|q \\ p\nmid m\endSb
\Big(1 - \frac1{(p-1)^2}\Big).
$$
Since $q$ is divisible by every prime $p\le h^2$, the above is
$$
= \prod_{p|m}\Big(1+\frac1{p-1}\Big)\prod_{p\nmid m}\Big(1 - \frac1{(p-1)^2}\Big)
+O\big(1/h^2\big).
$$
Here the main term is
$$
\prod_{p|m}\Big(1-\frac1p\Big)^{\!-2}\Big(1-\frac1p\Big)\times
\prod_{p\nmid m}\Big(1-\frac1p\Big)^{\!-2}\Big(1-\frac2p\Big)
= \eufm S(\{0,m\}).
$$
Goldston \cite{7} has shown that 
$$
2\sum_{m=1}^h (h-m)\eufm S(\{0,m\}) 
= h^2 - h\log h + Bh +O\big(h^{1/2+\varepsilon}\big),
\tag 47
$$
so we have the stated estimate.
\enddemo

    It is worth noting that (47) can also be written in the form
$$
\sum\Sb d_1, d_2 \\ 1\le d_i\le h \\ d_1 \ne d_2 \endSb
\eufm S(\{d_1,d_2\}) = h^2 - h\log h + Bh +O\big(h^{1/2+\varepsilon}\big).
\tag 48
$$      

    The term $d=1$ contributes $h^2 - h$ to the left hand side of (46).  Thus
if $q$ is divisible by every prime not exceeding $h^2$, then
$$
\sum \Sb d|q \\ d > 1\endSb 
\frac{\mu(d)^2}{\phi(d)^2}\Big(\sum\Sb a=1\\(a,d)=1\endSb^d
|E(a/d)|^2 \ - \phi(d)h\Big)\ = \ - h\log h + Ah+O\big(h^{1/2+\varepsilon}\big)
\tag 49
$$      
where $A = 2 - C_0 - \log 2\pi$.
\bigskip

     We now begin the main body of the proof of Theorem 2.  We apply (45) with
$y = h^{k+1}$, and set $Q = \prod_{p\le y}p$.  Thus
$$
R_k(h) = \sum\Sb q_1, \ldots, q_k \\ 1 < q_i \\ q_i|Q \endSb 
\ \prod_{i=1}^k\frac{\mu(q_i)}{\phi(q_i)}S(q_1,\ldots,q_k;h)\ +\ O(1)
\tag 50
$$       
where
$$
S(q_1,\ldots, q_k;h) = 
\sum \Sb d_1, \ldots, d_k \\ 1\le d_i \le h \\ d_i\text{ distinct} \endSb\
\sum\Sb a_1,\ldots, a_k \\ 1\le a_i \le q_i \\ 
(a_i,q_i)=1 \\ \sum a_i/q_i\in\Bbb Z\endSb
e\Big(\sum_{i=1}^k \frac{a_id_i}{q_i}\Big).
\tag 51
$$

     By comparing (50) with (12) we find that if the condition that the $d_i$ 
should be distinct were omitted, then the main term in (50) would be exactly 
$V_k(Q;h)$.  The bulk of our argument is devoted to an effort to remove this 
condition.  Put $\delta_{i,j} = 1$ if $d_i=d_j$, $\delta_{i,j}=0$ otherwise.  Thus
$$
\prod_{1\le i < j \le k}(1 - \delta_{i,j}) 
= \cases 1 &\text{if the $d_i$ are distinct;}\\
         0 &\text{otherwise.}
  \endcases
$$
When the left hand side above is expanded, we obtain a linear combination of products
of the $\delta$ symbols.  Let $\Delta$ denote such a product, and $|\Delta|$ the
number of factors in the product.  We define an equivalence relation on these 
$\delta$-products by setting $\Delta_1 \sim \Delta_2$ if $\Delta_1$ and $\Delta_2$
have the same value for all choices of $d_1, \ldots, d_k$.  For example,
$\delta_{1,2}\delta_{1,3}\sim\delta_{1,2}\delta_{2,3}
\sim\delta_{1,2}\delta_{1,3}\delta_{2,3}$.  Given a partition 
$\eusm P=\{\eusm S_1, \ldots, \eusm S_M\}$ of the set $\{1,\ldots, k\}$, let
$$
\Delta_{\eusm P} = \prod_{m=1}^M 
\prod\Sb i < j \\ i\in\eusm S_m \\ j \in\eusm S_m\endSb
\delta_{i,j}\,.
$$
We see easily that every equivalence class of $\delta$-products contains a unique
$\Delta_{\eusm P}$.  Thus we have a bijective correspondence between equivalence
classes of $\delta$-products and partitions of $\{1,\ldots, k\}$.  For a partition
$\eusm P$, put
$$
w(\eusm P) = \sum_{\Delta\sim\Delta_{\eusm P}}(-1)^{|\Delta|}.
$$
Thus
$$
\prod_{1\le i < j \le k}(1-\delta_{i,j}) = \sum_{\eusm P}w(\eusm P)\Delta_{\eusm P},
$$
and it follows that
$$
S(q_1,\ldots,q_k;h) = \sum_{\eusm P}w(\eusm P)
\sum\Sb a_1,\ldots,a_k\\ 1\le a_i\le q_i\\(a_i,q_i)=1\\ \sum a_i/q_i\in\Bbb Z\endSb
\prod_{m=1}^M E\Big(\sum_{i\in\eusm S_m}\frac{a_i}{q_i}\Big)
\tag 52
$$
where $\eusm P = \{\eusm S_1, \ldots, \eusm S_M\}$  

    If there is a prime $p$ such that $p|q_i$ for exactly one $i$, then the condition
$\sum a_i/q_i\in\Bbb Z$ cannot be fulfilled with $(a_i,q_i)=1$, so the sum (52) is
empty, and hence $S(q_1,\ldots, q_k; h) = 0$.  We therefore assume that each prime
dividing $[q_1,\ldots, q_k]$ divides at least two of the $q_i$.

    To facilitate our discussion of various types of partitions, we introduce some
notation.  Let $\eusm M = \{1, \ldots, M\}$.  For a partition 
$\eusm P = \{\eusm S_1,\ldots, \eusm S_M\}$, put
$$
\xalignat2
\eusm M_1 &= \{m\in\eusm M : \operatorname{card}\eusm S_m = 1\}, &m_1 
&=\operatorname{card} \eusm M_1; \\
\eusm M_2 &= \{m\in\eusm M : \operatorname{card}\eusm S_m \ge 2\}, &m_2
&=\operatorname{card}\eusm M_2; \\
\eusm N_1 &= \bigcup_{\text{card\,}\eusm S_m=1}\eusm S_m, &n_1 
&= \operatorname{card} \eusm N_1; \\
\eusm N_2 &= \bigcup_{\text{card\,}\eusm S_m\ge 2}\eusm S_m, &n_2
&= \operatorname{card} \eusm N_2.
\endxalignat 
$$              
Of course, $m_1+m_2 = M$, $n_1+n_2 = k$, and $m_1 = n_1$.  We first bound the
contribution made by those partitions such that $\operatorname{card}\eusm S_m\ge 3$
for some $m$.  For $m\in\eusm M_1$ we use (27) to see that 
$|E(a_i/q_i)|\le F(a_i/q_i)$.  For $m\in\eusm M_2$ we use the trivial bound
$$
\Big|E\Big(\sum_{i\in\eusm S_m}\frac{a_i}{q_i}\Big)\Big| \le h.
$$
Hence
$$
\Bigg|
\sum\Sb a_1,\ldots,a_k\\ 1\le a_i\le q_i\\(a_i,q_i)=1\\ \sum a_i/q_i\in\Bbb Z\endSb
\prod_{m=1}^M E\Big(\sum_{i\in\eusm S_m}\frac{a_i}{q_i}\Big)
\Bigg|
\ \le h^{m_2}
\sum\Sb a_1,\ldots,a_k\\ 1\le a_i\le q_i\\(a_i,q_i)=1\\ \sum a_i/q_i\in\Bbb Z\endSb
\ \prod_{i\in\eusm N_1}F\Big(\frac{a_i}{q_i}\Big).
$$
For $i\in \eusm N_1$ we take $G_i(x) = F(x)$ for $0 < x < 1$, and $G_i(1)=0$.  For
$i\in \eusm N_2$, we take $G_i(x)=1$ for all $x$.  Thus by Lemma~1 we see that 
the above is
$$
\le \frac{h^{m_2}}{[q_1,\ldots,q_k]}\prod_{i\in\eusm N_1}
\Big(q_i\sum_{a=1}^{q_i-1}F\Big(\frac a{q_i}\Big)^{\!2}\Big)^{\!1/2}
\prod_{i\in\eusm N_2}\Big(q_i\sum_{a=1}^{q_i}1\Big)^{\!1/2}.
\tag 53
$$
If $q \le h$, then
$$
\sum_{a=1}^{q-1}F\Big(\frac aq\Big)^{\!2} 
\ll \sum_{0<a\le q/2}\Big(\frac qa\Big)^{\!2}
\ll q^2.
$$
If $q > h$, then
$$
\sum_{a=1}^{q-1}F\Big(\frac aq\Big)^{\!2} \ll \sum_{0<a\le q/h}h^2\ 
+ \ \sum_{q/h<a\le q/2}\Big(\frac qa\Big)^{\!2} \ll qh\,.
$$
Thus in any case,
$$
\sum_{a=1}^{q-1}F\Big(\frac aq\Big)^{\!2} \ll q\min(q,h)\,.
\tag 54
$$
Thus the expression (53) is
$$
\ll \frac{q_1\cdots q_k}{[q_1,\ldots, q_k]}h^{n_1/2+m_2}.
$$
The $n_2$ members of $\eusm N_2$ are partitioned into $m_2$ sets, each one 
containing at least two members, and at least one containing $3$ or more numbers.
Thus
$$
n_2 = \sum_{m\in\eusm M_2}\operatorname{card}\eusm S_m 
\ge 1 + 2\sum_{m\in\eusm M_2}1 = 1 + 2m_2,
$$
and hence
$$
\frac{n_1}2 + m_2 \le \frac{n_1+n_2-1}2 = \frac{k-1}2.
$$
We also observe that
$$
\align
\sum\Sb q_1,\ldots, q_k \\ q_i|Q \endSb 
\bigg(\prod_{i=1}^k\frac{\mu(q_i)^2}{\phi(q_i)}\bigg)
\frac{q_1\cdots q_k}{[q_1,\ldots, q_k]}
&\le \sum_{d|Q}\frac1d\Big(\sum_{q|d}\frac q{\phi(q)}\Big)^{\!k}
= \sum_{d|Q}\prod_{p|d}\Big(1+\frac p{p-1}\Big)^{\!k} \\
&= \prod_{p|Q}\Big(1+\frac1p\Big(1+\frac p{p-1}\Big)^{\!k}\Big)
\ll_k (\log h)^{2^k}.
\endalign
$$
Thus we have shown that
$$
\aligned
R_k(h) &= \sum \Sb \eusm P \\ \operatorname{card}\eusm S_m\le 2\endSb
w(\eusm P)\sum\Sb q_1,\ldots, q_k \\ 1 < q_i|Q \endSb\ 
\prod_{i=1}^k\frac{\mu(q_i)}{\phi(q_i)}\sum
\Sb a_1,\ldots, a_k \\ 1\le a_i\le q_i \\ (a_i,q_i)=1\\ \sum a_i/q_i\in\Bbb Z\endSb
\prod_{m=1}^M E\Big(\sum_{i\in\eusm S_m}\frac{a_i}{q_i}\Big) \\ 
&\qquad +\ O\big(h^{(k-1)/2+\varepsilon}\big).
\endaligned
\tag 55
$$

    Suppose that the partition $\eusm P$ consists of $j$ doubleton sets and $k-2j$
singleton sets.  Since no other $\delta$-product is equivalent to 
$\Delta_{\eusm P}$, and $|\Delta_{\eusm P}|=j$, so $w(\eusm P) = (-1)^j$.  The 
number of such partitions is 
$$
\Big({k\atop 2j}\Big)\frac{(2j)!}{j!2^j}.
$$      
Since the $q_i$ are interchangeable, we multiply by the above factor, and restrict 
our attention to one such partition: doubletons $\{i, i+j\}$ for $1\le i\le j$ 
and singletons
$\{i\}$ for $2j+1\le i\le k$.  Thus the main term in (55) is
$$
\sum_{0\le j\le k/2}(-1)^j\Big({k\atop 2j}\Big)\frac{(2j)!}{j!2^j}
\sum\Sb q_1,\ldots, q_k \\ 1 < q_i|Q \endSb\ 
\prod_{i=1}^k\frac{\mu(q_i)}{\phi(q_i)}\sum
\Sb a_1,\ldots, a_k \\ 1\le a_i\le q_i \\ (a_i,q_i)=1\\ \sum a_i/q_i\in\Bbb Z\endSb
\prod_{i=1}^jE\Big(\frac{a_i}{q_i}+\frac{a_{i+j}}{q_{i+j}}\Big)
\prod_{i=2j+1}^k \!\!E\Big(\frac{a_i}{q_i}\Big)\,.
\tag 56
$$

    For $1\le i\le j$, let $b_i$ and $r_i$ be defined by the relations
$$
\frac{b_i}{r_i} \equiv \frac{a_i}{q_i}+ \frac{a_{i+j}}{q_{i+j}}\pmod 1,\qquad
1\le b_i\le r_i, \qquad (b_i, r_i) = 1,
$$      
and put
$$
H\Big(\frac br\Big) = E\Big(\frac br\Big)\sum\Sb d_1, d_2 \\ 1<d_i|Q \endSb
\frac{\mu(d_1)\mu(d_2)}{\phi(d_1)\phi(d_2)}
\sum\Sb c_1,c_2 \\ 1\le c_i \le d_i \\ (c_i,d_i)=1 \\ 
\frac{c_1}{d_1}+\frac{c_2}{d_2}\equiv \frac br \pmod 1\endSb 1\,.
\tag 57
$$
Then the sum over the $q_i$ in (56) is 
$$
\sum\Sb r_1,\ldots,r_j \\ r_i|Q \endSb 
\sum \Sb b_1,\ldots, b_j \\ 1\le b_i\le r_i \\ (b_i,r_i)=1\endSb 
\ \prod_{i=1}^jH\Big(\frac{b_i}{r_i}\Big)
\sum\Sb q_{2j+1},\ldots,q_k \\ 1 < q_i|Q \endSb \
\sum\Sb a_{2j+1},\ldots, a_k \\ 1\le a_i \le q_i \\ (a_i,q_i)=1 \\ 
\sum a_i/q_i\in\Bbb Z\endSb \ \prod_{i=2j+1}^k \frac{\mu(q_i)}{\phi(q_i)}
E\Big(\frac{a_i}{q_i}\Big)\,.
$$
We now separate those $i$ for which $r_i = 1$ from those with $r_i > 1$.  Let $\ell$
denote the number of $i$ for which $r_i > 1$.  Since there are 
$\big({j\atop \ell}\big)$ ways of choosing the $\ell$ values of the $i$ from 
$\{1, \ldots, j\}$, the above is
$$
\sum_{\ell=0}^j \Big({j\atop \ell}\Big)H(1)^{j-\ell}M(\ell)
\tag 58
$$
where
$$
M(\ell) = \sum\Sb r_1,\ldots, r_\ell \\ 1<r_i|Q \endSb 
\sum \Sb b_1,\ldots, b_\ell \\ 1\le b_i\le r_i \\ (b_i,r_i)=1\endSb
\prod_{i=1}^\ell H\Big(\frac{b_i}{r_i}\Big)
\sum \Sb q_{2j+1},\ldots, q_k \\ 1<q_i|Q\endSb
\sum \Sb a_{2j+1},\ldots, a_k \\ 1\le a_i \le q_i \\ (a_i,q_i)=1\\
\sum a_i/q_i + \sum b_i/r_i \in\Bbb Z \endSb 
\prod_{i=2j+1}^k\frac{\mu(q_i)}{\phi(q_i)}E\Big(\frac{a_i}{q_i}\Big).
$$
We note that $M(0) = V_{2k-j}(Q;h)$.  

    Next we show that the contributions of $\ell > 0$ can be absorbed in the error 
term.  If there is a prime $p$ that divides exactly one of the numbers 
$r_1, \ldots, r_\ell,q_{2j+1},\ldots, q_k$, then the condition that 
$\sum a_i/q_i+\sum b_i/r_i\in\Bbb Z$ cannot be fulfilled with 
$(a_i,q_i)=(b_i,r_i)=1$, so the sum is empty, and the sum over the $q_i$ and $r_i$ 
vanishes.  Thus we may restrict our attention to those choices of $q_i$ and $r_i$ 
for which every prime divisor of $d = [q_{2j+1},\ldots, q_k,r_1,\ldots, r_\ell]$ 
divides at least two of these numbers.  Hence by Lemma~1, 
$$
M(\ell) \ll \sum\Sb r_1, \ldots, r_\ell \\ 1<r_i|Q \endSb\,
\sum\Sb q_{2j+1},\ldots, q_k \\ 1 < q_i|Q \endSb \frac1d\prod_{i=1}^\ell\!
\Big(r_i\sum_{b=1}^{r_i-1}\Big|H\Big(\frac b{r_i}\Big)\Big|^2\Big)^{\!1/2}\!\!
\prod_{i=2j+1}^k\!\!\Big(\frac {q_i}{\phi(q_i)^2}\sum_{a=1}^{q_i-1}
F\Big(\frac a{q_i}\Big)^{\!2}
\Big)^{\!1/2}.
\tag 59
$$
          
          In order to assess the above, we estimate $H(b/r)$.  We note that $H(b/r)=0$
if $r\nmid Q$.  Thus we suppose that $r|Q$, and that $r > 1$.  For $i = 1, 2$ we 
write $d_i = s_it_i$ where $s_i|r$, $t_i|Q/r$.  By the Chinese Remainder Theorem 
there exist unique $e_i$ (mod $s_i$) and $f_i$ (mod $t_i$) such that
$$
\frac{e_i}{s_i} + \frac{f_i}{t_i} \equiv \frac{c_i}{d_i}\pmod 1,\qquad
(e_i,s_i) = (f_i,t_i) = 1.  
$$        
From the conditions $c_1/d_1 + c_2/d_2\equiv b/r \pmod 1$, $(b,r)=1$ it follows that
$$
\frac{e_1}{s_1}+\frac{e_2}{s_2}\equiv \frac br\pmod 1, \qquad [s_1,s_2]=r, \qquad
t_1 = t_2, \qquad f_1+f_2\equiv 0 \pmod{t_1}.
$$
Put $t = t_1=t_2$.  Hence
$$
H\Big(\frac br\Big)=E\Big(\frac br\Big)\sum\Sb s_1, s_2 \\ [s_1,s_2]=r\endSb
\sum \Sb e_1,e_2 \\ 1\le e_i\le s_i \\ (e_i,s_i)=1 \\ 
\frac{e_1}{s_1}+\frac{e_2}{s_2}\equiv \frac br\pmod 1\endSb
\frac{\mu(s_1)\mu(s_2)}{\phi(s_1)\phi(s_2)}\sum_{t|Q/r}\frac{\mu(t)^2}{\phi(t)}\,.
$$
For given $s_1, s_2$, the number of pairs $e_1, e_2$ with the required properties is
$\le \phi((s_1,s_2))= \phi(s_1)\phi(s_2)/\phi(r)$.  Hence
$$
\align
H\Big(\frac br\Big) &\ll F\Big(\frac br\Big)\sum\Sb s_1, s_2 \\ [s_1,s_2]=r\endSb
\frac1{\phi(r)}\prod_{p|Q/r}\Big(1 + \frac1{p-1}\Big)
= F\Big(\frac br\Big)\frac{3^{\omega(r)}}r\prod_{p|Q}\Big(1+\frac1{p-1}\Big) \\
&\ll F\Big(\frac br\Big)\frac{3^{\omega(r)}}r\log h\,.
\endalign
$$
By (54) it follows that
$$
\sum_{b=1}^{r-1}\Big|H\Big(\frac br\Big)\Big|^2 
\ll \frac hr 9^{\omega(r)}(\log h)^2\,.
$$
On inserting this and (54) in (59), we find that
$$
\align
M(\ell) &\ll \sum\Sb r_1,\ldots, r_\ell \\ r_i|Q \endSb 
\sum \Sb q_{2j+1}\ldots, q_k \\ q_i|Q \endSb \frac 1d
\prod_{i=1}^\ell\big(h^{1/2}3^{\omega(r_i)}\log h\big)
\prod_{i=2j+1}^k\frac {h^{1/2}q_i}{\phi(q_i)} \\
&\ll h^{(k-2j+\ell)/2}(\log h)^{\ell}\sum_{d|Q}\frac1d
\Big(\sum_{r|d}3^{\omega(r)}\Big)^{\!\ell}
\Big(\sum_{q|d}\frac q{\phi(q)}\Big)^{\!k-2j} \\
&\le h^{(k-2j+\ell)/2}(\log h)^{\ell}
\prod_{p|Q}\Big(1+\frac{4^{k-2j+\ell}}p\Big) \\
&\ll h^{(k-2j+\ell)/2+\varepsilon}.
\endalign
$$
From (57) we see that
$$
H(1) = h\sum \Sb d|Q \\ d > 1 \endSb \frac{\mu(q)^2}{\phi(q)}.
$$
Hence the expression (58) is
$$
\Big(h\sum\Sb d|Q \\ d > 1 \endSb \frac{\mu(q)^2}{\phi(q)}\Big)^{\!j}V_{k-2j}(Q;h)
\ + \ O\big(h^{(k-1)/2+\varepsilon}\big)\,.
$$
On inserting this in (56), we find that
$$
R_k(h) = \sum_{0\le j\le k/2} \Big({k\atop 2j}\Big)\frac{(2j)!}{j!2^j}
\Big(\!-h\sum\Sb d|Q \\ d > 1\endSb \frac{\mu(d)^2}{\phi(d)}\Big)^{\!j}V_{k-2j}(Q;h)
\ + \ \big(h^{(k-1)/2+\varepsilon}\big)\,.
\tag 60
$$

    We are at last prepared to make our appeal to Theorem 1.  If $k$ is odd then 
$k-2j$ is odd, so there is no main term.  Suppose that $k$ is even.  Then the 
main term is
$$
\align
&\ \sum_{0\le j\le k/2}\Big({k\atop 2j}\Big)\frac{(2j)!}{j!2^j}
\Big(\!-h\sum\Sb d|Q \\ d>1\endSb \frac{\mu(d)^2}{\phi(d)}\Big)^{\!j}
\frac{(k-2j)!}{(k/2-j)!2^{k/2-j}}V_2(Q;h)^{k/2-j} \\
&=\frac{k!}{(k/2)!2^{k/2}}\sum_{j=0}^{k/2}\Big({k/2\atop j}\Big)V_2(Q;h)^{k/2-j}
\Big(\!-h\sum\Sb d|Q \\ d > 1\endSb \frac{\mu(d)^2}{\phi(d)}\Big)^{\!j} \\
&= \mu_k\Big(V_2(Q;h) - h\sum\Sb d|Q \\ d>1\endSb
\frac{\mu(d)^2}{\phi(d)}\Big)^{\!k/2}
\endalign
$$
by the binomial theorem.

     By (49) we see that the above is 
$$
\mu_k(-h\log h + Ah)^{k/2} + O\big(h^{(k-1)/2+\varepsilon}\big)\,.
$$       
This gives the stated result.

\head
3. Proof of Theorem 3
\endhead

\noindent
Clearly $M_0(N;H)=N$.  Thus we have the case $K=0$ of (21) unconditionally, and
with no error term.  It is also convenient to dispose of the case $K=1$ before
proceeding to the main argument.  Since $\eufm S(\{h\})=1$, it follows from our
hypothesis (20) that 
$\sum_{n=1}^N \Lambda(n+h) = N +O\big(N^{1/2+\varepsilon}\big)$.
Hence $\sum_{n=1}^N \Lambda_0(n+h)\ll N^{1/2+\varepsilon}$, and thus by (19)
we see that $M_1(N;H)\ll HN^{1/2+\varepsilon}$, which suffices.  

     From now on we assume that $K$ is fixed, $K\ge 2$, and we ignore possible
dependence of implicit contstants on $K$.  Let $\eusm D = \{d_1, \ldots, d_k\}$
be a set of $k$ distinct integers with $1\le d_i\le H$ for $1\le i\le k$.  Suppose
that the $h_i$ in (19) take the values $d_i$ with multiplicities $M_i$.  Then
the right hand side of (19) is 
$$
= \ \sum_{k=1}^K \ \sum\Sb M_1,\ldots, M_k \\ M_i\ge 1\\\sum M_i = K\endSb
\Big({K\atop M_1\ \cdots\ M_k}\Big)\frac1{k!}
\sum\Sb d_1,\ldots, d_k \\ 1\le d_i\le H \\ d_i\text{ distinct}\endSb \sum_{n=1}^N \
\prod_{i=1}^k \Lambda_0(n+d_i)^{M_i}.
\tag 61
$$
Here the $1/k!$ is necessary because any permutation of $d_1, \ldots, d_k$ gives rise
to the same set $\eusm D$.  For positive integers $m$ we put 
$\Lambda_m(n) = \Lambda(n)^m\Lambda_0(n)$.  If $M\ge 1$, then by the binomial theorem
$$
\Lambda_0(n)^M = \Lambda_0(n)(\Lambda(n)-1)^{M-1}
= \sum_{m=0}^{M-1}(-1)^{M-m-1}\Big({M-1\atop m}\Big)\Lambda_m(n).
$$
On inserting this in (61), we find that
$$
\aligned
M_K(N;H) &= \sum_{k=1}^K\frac1{k!}
\sum\Sb M_1,\ldots, M_k\\ M_i\ge 1\\ \sum M_i=K\endSb
\Big({K\atop M_1\ \cdots M_k}\Big) \\
&\qquad \times\sum\Sb m_1,\ldots, m_k \\ 0\le m_i < M_i\endSb
\prod_{i=1}^k (-1)^{M_i-1-m_i}\Big({M_i-1\atop m_i}\Big)L_k(\boldkey m)
\endaligned
\tag 62
$$
where
$$
L_k(\boldkey m) = 
\sum\Sb d_1,\ldots, d_k \\ 1\le d_i\le H \\ d_i\text{ distinct}\endSb
\sum_{n=1}^N \ \prod_{i=1}^k \Lambda_{m_i}(n+d_i).
\tag 63
$$
To estimate the $L_k(\boldkey m)$, we must distinguish between those $i$ for which
$m_i =0$ and those for which $m_i > 0$.  To this end we set 
$\eusm K = \{1,\ldots k\}$, and introduce the following notation:
$$
\xalignat2
\eusm H &= \{i\in\eusm K : m_i \ge 1\}, & h &= \operatorname{card}\eusm H; \\
\eusm I &= \{i\in\eusm K : m_i = 0\}, & k-h &= \operatorname{card}\eusm I; \\
\eusm J &\subseteq \eusm K, & j &= \operatorname{card} \eusm J.
\endxalignat
$$
Thus
$$
\align
\prod_{i\in\eusm I}\Lambda_0(n+d_i)\prod_{i\in\eusm H}\Lambda(n+d_i)
&= \prod_{i\in\eusm I}\Lambda_0(n+d_i)\prod_{i\in\eusm H}(\Lambda_0(n+d_i)+1) \\
&= \sum\Sb \eusm J \\ \eusm I\subseteq \eusm J\subseteq \eusm K\endSb \
\prod_{i\in \eusm J}\Lambda_0(n+d_i).
\endalign
$$
From our hypothesis (20) it follows that
$$
\sum_{n\le x}\ \prod_{i\in\eusm I}\Lambda_0(n+d_i)\prod_{i\in \eusm H}\Lambda(n+d_i)
= x\sum\Sb \eusm J \\ \eusm I \subseteq \eusm J \subseteq \eusm K \endSb
\eufm S_0(\eusm D_{\eusm J})
+O\big(N^{1/2+\varepsilon}\big)
$$
uniformly for $0\le x\le N$ where $\eusm D_{\eusm J} = \{d_i:i\in\eusm J\}$.
With $\boldkey m = (m_1, \ldots, m_k)$ fixed for the
moment, write the above briefly as $f(x)= cx +O\big(N^{1/2+\varepsilon}\big)$.
Then
$$
\align
\sum_{n=1}^N &\Big(\prod_{i\in\eusm I}\Lambda_0(n+d_i)\Big)
\Big(\prod_{i\in\eusm H}\Lambda(n+d_i)(\log(n+d_i))^{m_i-1}
(\log(n+d_i) - 1)\Big)
\tag 64 \\
&= \int_{1^-}^N \prod_{i\in\eusm H}(\log (x+d_i))^{m_i-1}
(\log (x+d_i) - 1)\,df(x), \\
\intertext{which by integration by parts is}
&= c\int_1^N\prod_{i\in\eusm H}(\log(x+d_i))^{m_i-1}(\log(x+d_i)- 1)\,dx 
+O\big(N^{1/2+\varepsilon}\big).
\tag 65 \\
\endalign
$$
For $m > 0$, $\Lambda_m(n)$ is nonzero only when $n$ is a primepower, and 
$$
\Lambda_m(n)= \Lambda(n)(\log n)^{m-1}(\log n - 1)
$$
if $n$ is prime.  Thus if
$n$ is an integer such that $n+d_i$ is prime for all $i\in\eusm H$, then the summand
in (64) is
$$
\prod_{i=1}^k \Lambda_{m_i}(n+d_i).
$$
Those $n$ for which $n+d_i$ is a higher power of a prime for one or more 
$i\in\eusm H$ contribute an amount $\ll N^{1/2+\varepsilon}$ to the sum (64).  
Thus the sum (64) is
$$
= \sum_{n=1}^N\ \prod_{i=1}^k\Lambda_{m_i}(n+d_i)\ 
+ \ O\big(N^{1/2+\varepsilon}\big).
\tag 66
$$
Next we approximate the integral in (65) by a similar integral that is 
independent of the $d_i$.  First we note that if $k = K$, then $M_i = 1$ for
all $i$,  and hence $m_i=0$ for all $i$, so that $h = 0$.  Thus we may suppose
that $k < K$.  It is useful to note that
$$
\sum_{i\in\eusm H} m_i = \sum_{i=1}^k m_i \le \sum_{i=1}^k (M_i-1) = K-k\,.
\tag 67
$$   
If $x\ge 1$ and $1\le d \le H$, then
$$
\log(x+d) = \log x + O(d/x) = \log x +O(H/x).
$$
Thus the integrand is
$$
\prod_{i\in\eusm H}\big((\log x)^{m_i-1}(\log x - 1)\big) 
+O\big(Hx^{-1}(\log N)^{K-k-1}\big),
$$
and so the integral is $I_{\boldkey m}(N)+O\big(H(\log N)^{K-k}\big)$ where
$$
I_{\boldkey m}(N) = \int_1^N \prod_{i\in\eusm H}
\big((\log x)^{m_i-1}(\log x - 1)\big)\,dx.
$$
On assembling our estimates, we find that
$$
\sum_{n=1}^N \ \prod_{i=1}^k\Lambda_{m_i}(n+d_i) \
= \ \bigg(\sum\Sb \eusm J \\ \eusm I \subseteq \eusm J \subseteq \eusm K\endSb
\eufm S_0(\eusm D_{\eusm J})\bigg)\big(I_{\boldkey m}(N)
+O\big(H(\log N)^{K-k}\big)\big)
+O\big(N^{1/2+\varepsilon}\big).
$$
By (43) we see that the first error term is
$\ll H(\log N)^K \ll N^{1/K+\varepsilon}\ll N^{1/2+\varepsilon}$, since 
$H\le N^{1/K}$ and $K\ge 2$.

     On summing both sides of the above over all choices of distinct $d_i$, we 
find that
$$
\align
L_k(\boldkey m) &= I_{\boldkey m}(N)
\sum\Sb \eusm J \\ \eusm I\subseteq \eusm J\subseteq \eusm K\endSb
\sum\Sb d_1, \ldots, d_k\\ 1\le d_i \le H \\ d_i \text{ distinct}\endSb
\eufm S_0(\eusm D_{\eusm J}) + O\big(H^kN^{1/2+\,\varepsilon}\big).
\intertext{Once the $d_i$ have been chosen for $i\in\eusm J$, there are
$(H-j)(H-j-1)\cdots(H-k+1)$ ways of choosing the remaining $d_i$.  Hence the 
above is}
&=I_{\boldkey m}(N)
\sum\Sb \eusm J \\ \eusm I\subseteq \eusm J\subseteq \eusm K\endSb
R_j(H)(H-j)\cdots(H-k+1) \ + \ O\big(H^kN^{1/2+\varepsilon}\big).
\endalign
$$       
The product in the sum is $H^{k-j}+O\big(H^{k-j-1}\big)$.  By Theorem~2 we know
that $R_j(H)\ll (H\log H)^{j/2}$, and using this we see that (isolating the term  
$\eusm J = \eusm I$) 
$$
\align 
\sum\Sb \eusm J \\ \eusm I\subseteq \eusm J\subseteq \eusm K\endSb
R_j(H)(H-j)&\cdots(H-k+1) 
\\
&= R_{k-h}(H)(H^h + O(H^{h-1})) + O\Big( \sum_{k-h+1 \le 
j\le k } (H\log H)^{j/2} H^{k-j}\Big)\\
&= R_{k-h}(H) H^h +O\big( (H\log H)^{(k-h+1)/2} H^{h-1}\big).\\
\endalign
$$
Therefore 
$$
L_k(\boldkey m) = I_{\boldkey m}(N)\Big(R_{k-h}(H)H^h 
+O\big((H\log H)^{(k-h+1)/2}H^{h-1}\big)\Big)
+O\big(H^kN^{1/2+\varepsilon}\big).
\tag 68
$$

We insert the above in (62).  In assessing the sizes of the various terms, it is
useful to note that
$$
K = \sum_{i=1}^k M_i = \sum_{i\in\eusm H}M_i \ + \ \sum_{i\in\eusm I}M_i 
  \ge 2\operatorname{card} \eusm H + \operatorname{card} \eusm I = 2h + (k-h) = h+k.
\tag 69  
$$
By (67) we see that 
$$
I_{\boldkey m}(N) \sim N (\log N)^{\sum_{i\in \eusm H} m_i } \ll N (\log N)^{K-k}. 
\tag 70
$$

First we show that terms for which $h+k < K$ contribute a negligible amount to (62).
Since $R_{k-h}(H) \ll (H\log H)^{(k-h)/2}$ by Theorem~2, it follows from
(68), (69), and (70) that
$$
\align
L_k(\boldkey m) &\ll N(\log N)^{K-k} (H\log H)^{(k-h)/2} H^{h}
 + H^kN^{1/2+\varepsilon}\\
&\ll N(\log N)^{K} \Big(\frac{H}{\log N}\Big)^{\!(k+h)/2} 
\Big(\frac{\log H}{\log N}\Big)^{\!(k-h)/2}
 + H^KN^{1/2+\varepsilon}.\\
\endalign
$$
Thus the contribution of these terms to (62) is 
$$
\ll N (\log N)^{K} \Big(\frac{H}{\log N}\Big)^{\!(K-1)/2} 
+ H^K N^{1/2+\varepsilon}. 
\tag 71
$$

     Finally we consider those terms in (62) for which $h+k = K$.  Since $h\le k$,
it follows that $k\ge K/2$.  We also have $\operatorname{card} \eusm H = h = K-k$, 
and $\operatorname{card} \eusm I = k-h = 2k - K$.  Since equality
holds in (69), it follows that $M_i = 2$ for all $i\in\eusm H$ and that $M_i = 1$
for all $i\in \eusm I$.  Thus $m_i = M_i - 1$ for all $i$, and for such $\boldkey m$
we have
$$
I_{\boldkey m}(N) = \int_1^N (\log x - 1)^{K-k}\,dx = I_{K-k}(N),
$$       
say.  Hence 
$$
\align
L_k(\boldkey m) &= 
I_{K-k}(N)R_{2k-K}(H)H^{K-k} 
+O\big(N(\log N)^{K-k} (H\log H)^{(2k-K+1)/2}H^{K-k-1}\big)
\\
&\hskip 1.76 in +O(H^k N^{1/2+\varepsilon})\\ 
&= I_{K-k}(N)R_{2k-K}(H)H^{K-k} +O\Big( N(\log N)^K 
\Big(\frac{H}{\log N}\Big)^{\!(K-1)/2}\Big) 
+ O(H^K N^{1/2+\varepsilon}).\\
\endalign
$$
Once $k$ is selected, there are precisely $\big({k\atop K-k}\big)$ ways of choosing
the set $\eusm H$, and hence using (71) and the above,
$$
\aligned
M_K(N;H) &= \sum_{K/2\le k\le K} \frac{K!}{k!2^{K-k}}\Big({k\atop K-k}\Big)
I_{K-k}(N)R_{2k-K}(H)H^{K-k} \\
&\qquad+O\Big( N(\log N)^K \Big(\frac{H}{\log N}\Big)^{\!(K-1)/2}\Big) 
 + O\big(H^{K}N^{1/2+\varepsilon}\big).
\endaligned
\tag 72
$$

     Suppose that $K$ is odd.  Then so also is $2k-K$, and hence by Theorem~2 the
main terms in (72) are 
$$
\ll N(\log N)^{K/2} H^{K/2-1/(7K)+\varepsilon}.  
$$
Thus in this case 
$$
\align
M_K(N;H) &\ll N(\log N)^{K/2} H^{K/2-1/(7K)+\varepsilon} +  
N(\log N)^K \Big(\frac{H}{\log N}\Big)^{\!(K-1)/2} + H^K N^{1/2+\varepsilon}\\
&\ll N(\log N)^{K/2}H^{K/2} \Big(\frac{H}{\log N}\Big)^{\!-1/(8K)} 
+H^K N^{1/2+\varepsilon}.
\\
\endalign
$$

     Suppose that $K$ is even.  By Theorem~2 it follows that
$$
\align
M_K(N;H) &= H^{K/2}\sum_{k=K/2}^K \frac{K!}{k!2^{K-k}}\binom{k}{K-k}
\mu_{2k-K}I_{K-k}(N)(-\log H+A)^{k-K/2} \\ 
&\qquad+O\Big( N(\log N)^{K/2} H^{K/2} \Big(\frac{H}{\log N}\Big)^{\!-1/(8K)}
 + H^K N^{1/2+\varepsilon}\Big).\\
\endalign
$$       
Since $\mu_k = k!/((k/2)!2^{k/2})$ when $k$ is even, on writing $k = K/2 + \ell$
we see that the main term above is 
$$
\mu_KH^{K/2}\sum_{\ell=0}^{K/2}\Big({K/2\atop \ell}\Big)I_{K/2-\ell}(N)
(-\log H+A)^{\ell}.
$$
On taking the sum inside the integral, we obtain the stated result by the
binomial theorem.

\head 
4. Heuristics in the manner of Hardy \& Littlewood
\endhead
  
\noindent
The reasoning here is conducted in the manner that Hardy
\& Littlewood \cite{10} used to formulate their quantitative version
of the prime $k$-tuple hypothesis.  Let
$$
S(\alpha) = \sum_{n=1}^N \Lambda(n)e(n\alpha).
$$
By the prime number theorem for arithmetic progressions we know that
if $(a,q)=1$, then $S(a/q) \sim \mu(q)N/\phi(q)$ provided that $q$ is 
not too 
large as a function of $N$.  By partial summation it follows that
$S(\alpha)\sim \mu(q)M(\alpha-a/q)/\phi(q)$ for $\alpha$ near $a/q$,
where $M(\beta) = \sum_{n=1}^N e(n\beta)$.  Put 
$E(\alpha)=\sum_{m=1}^h e(m\alpha)$, as in the proof of Theorem~1.
If $h$ is small compared with $N$, then $S(\alpha)E(-\alpha)$ is 
approximately
$$
\sum_{n=1}^N\Big(\sum_{m=1}^h \Lambda(m+n)\Big)e(n\alpha)\,.
$$
Moreover, when $\alpha$ is near $a/q$, $E(-\alpha)$ is approximately
$E(-a/q)$.  Thus when $\alpha$ is a number for which the expression
above is large, we expect that that it is approximately
$$
\sum_{1\le q\le N}\frac{\mu(q)}{\phi(q)}\sum \Sb a=1\\(a,q)=1\endSb^q
E(-a/q)M(\alpha-a/q)\,.
$$
We subtract the 
contribution of the term $a = q = 1$ from both sides to see that
$$
\sum_{n=1}^N\Big(\sum_{m=1}^h\Lambda(m+n) - h\Big)e(n\alpha)
\doteqdot\sum_{1< q\le N}\frac{\mu(q)}{\phi(q)}\sum \Sb a=1\\(a,q)=1\endSb^q
E(a/q)M(\alpha-a/q).
$$
Let $F(\alpha)$ denote the left hand side above.  The $k$-fold convolution
of $F$ with itself is
$$
\sum_{n=1}^N \Big(\sum_{m=1}^h \Lambda(m+n)\ - h\Big)^{\!k}e(n\alpha)
= \mathop{\int\cdots \int}_{\sum\alpha_i = \alpha\,(\text{mod}\,1)}
\prod_{i=1}^k F(\alpha_i)\,d\alpha_1\cdots d\alpha_k\,.
$$
We set $\alpha = 0$, and follow Hardy
\& Littlewood in assuming that the main term arises by the alignment of
the peaks in the multiple integral on the right hand side.  Thus we expect 
that
$$
\sum_{n=1}^N\Big(\sum_{m=1}^h \Lambda(m+n) - h\Big)^{\!k}
\sim N\sum\Sb q_1,\ldots, q_k\\ 1 < q_i \le N\endSb
\Big(\prod_{i=1}^k \frac{\mu(q_i)}{\phi(q_i)}\Big)
\sum \Sb a_1, \ldots, a_k \\ 1\le a_i\le q_i \\ (a_i,q_i)=1 \\ 
\sum a_i/q_i\in \Bbb Z\endSb\
\prod_{i=1}^k E(a_i/q_i).
$$
We note that $|E(\beta)|\asymp h$ if $\|\beta\|\le 1/h$, and that
$E(\beta)\ll 1/\|\beta\|$ if $\|\beta\|\ge 1/h$.  The asymptotic size of
the right hand side above could be determined by using the techniques used
to prove Theorem~1.  At this point we are content to argue more informally.
If $k$ is odd then the terms  do not
make a very significant contribution.  On the other hand, when $k$ is even, 
we find `diagonal terms' in which the $q_i$
are equal in pairs, with the corresponding $a_i$ being the negatives of each
other.  The pairings can be made in
$$
(k-1)(k-3)\cdots3\cdot1 = \mu_k
$$
ways, so the contribution of these terms is
$$
\mu_kN\!\sum\Sb q_1,\ldots,q_{k\!/2}\\ 1 < q_i \le N\endSb\,\prod_{i=1}^{k/2}\!\!
\bigg(\frac{\mu(q_i)^2}{\phi(q_i)^2}\!\!\sum\Sb a_i=1\\(a_i,q_i)=1\endSb^{q_i}\!\!
|E(a_i/q_i)|^2\!\bigg)\
= \ \mu_kN\bigg(\sum_{1 < q \le N}\frac{\mu(q)^2}{\phi(q)^2}\!
\sum\Sb a=1\\(a,q)=1\endSb^q\! |E(a/q)|^2\!\bigg)^{\!k/2}.
$$
If there are further equalities among the $q_i$ beyond this pairing, then the
combinatorics must be adjusted, but such configurations contribute
a lesser amount.  Likewise, that the non-diagonal terms contribute a lesser 
amount 
can be demonstrated by using the techniques we used to prove Theorem~1. 
If $q < h$ then the inner sum above is 
$\ll \sum_{0<a<q}\|a/q\|^{-2}\ll q^2$, 
but if $q \ge h$ then the inner
sum is approximately $h\phi(q)$.  Since
$$
\sum_{q\le y} \frac{\mu(q)^2}{\phi(q)} = \log y +O(1),
$$
it follows that the expression to be estimated is
$$
\sim \mu_k N\Big(h\log \frac Nh\Big)^{\!k/2},
$$
which supports the Conjecture.

\bigskip

\line{\smc Dept.\ of Math., Univ.\ of Michigan, Ann Arbor, MI 48109--1109, USA
\hfill}
\line{email: {\tt hlm\@umich.edu}\hfill}
\medskip
\line{\smc Dept.\ of Math., Univ.\ of Michigan, Ann Arbor, MI 48109--1109, USA\hfill}
\line{email: {\tt ksound\@umich.edu}\hfill}

\enddocument